\documentclass[12pt]{article}
\usepackage{amsmath, amssymb, amsthm,amscd}

\newcommand{\1}{\mathbf {1}}
\newcommand{\Z}{{\mathbb Z}}

\newcommand{\C}{{\mathbb C}}

\newcommand{\h}{{\mathfrak h}}

\newcommand{\mmbox}[1]{\mbox{\scriptsize #1}}
\newtheorem{definition}{Definition}
\newtheorem{lemma}{Lemma}
\newtheorem{theorem}{Theorem}
\newtheorem{corollary}[theorem]{Corollary}
\title{On intertwining operators and 
finite automorphism groups of vertex operator algebras}
\author{Kenichiro Tanabe\footnote{Partially supported by JSPS
Grant-in-Aid for Scientific Research No. 14740061.}\\\\
Institute of Mathematics \\
University of Tsukuba
}
\date{}
\begin{document}
\maketitle
\begin{abstract}
Let $V$ be a simple vertex operator algebra and 
$G$ a finite automorphism group.  
We give a construction of intertwining operators
for irreducible $V^G$-modules which occur as
submodules of irreducible $V$-modules by using  
intertwining operators for $V$.

We also determine some fusion rules for a vertex operator algebra
as an application.
\end{abstract}

\section{Introduction}

Let $V$ be a simple vertex operator algebra (cf. \cite{B}, \cite{FHL} and \cite{FLM}), and 
$G$ a finite automorphism group. 
It is an important problem to understand the module category 
for the vertex operator algebra $V^G$ of $G$-invariants.
In \cite{DVVV}, this question was asked and 
several ideas are proposed.

For a simple vertex operator algebra $V$, it is shown in \cite{DY} 
that every irreducible $V$-module is 
a completely reducible $V^G$-module as a natural consequence
of a duality theorem of Schur-Weyl type.
In this paper
we give a construction of intertwining operators
for irreducible $V^G$-modules which occur as submodules 
of irreducible $V$-modules by using  intertwining operators for $V$.

Let's state our results more explicitly. 
Firstly, we need to recall the results of Dong--Yamskulna\cite{DY}.
For 
an irreducible $V$-module $(L,Y_L)$ and $a\in G$
we define a new irreducible $V$-module $(L\circ a,Y_{L\circ a})$.
Here $L\circ a$ is equal to $L$ and $Y_{L\circ a}(u,z)=Y_L(au,z)$.
Let ${\cal S}$ 
be a finite set of inequivalent irreducible 
$V$-modules which is closed under this right action of $G$. 
In \cite{DY} they define a finite dimensional semisimple associative algebra 
${\cal A}_{\alpha}(G,{\cal S})$ over ${\mathbb C}$
and show a duality theorem of Schur-Weyl type
for the actions 
of $V^G$ and ${\cal A}_{\alpha}(G,{\cal S})$ on 
the direct sum of $V$-modules
in ${\cal S}$ which is denoted by ${\cal L}$.
That is, as a ${\cal A}_{\alpha}(G,{\cal S})\otimes V^G$-module,
${\cal L}=\bigoplus_{(j,\lambda)\in \Gamma}
W_{(j,\lambda)}\otimes M_{(j,\lambda)}$,
where $\{W_{(j,\lambda)}\}_{(j,\lambda)\in\Gamma}$ is the set of
all inequivalent irreducible ${\cal A}_{\alpha}(G,{\cal S})$-modules and 
$M_{(j,\lambda)}$ is the multiplicity spaces of 
$W_{(j,\lambda)}$ in ${\cal L}$.
Each $M_{(j,\lambda)}$ is a nonzero irreducible $V^G$-module
and the different multiplicity spaces are inequivalent 
$V^G$-modules.
In this paper we consider intertwining operators 
for irreducible $V^G$-modules constructed from irreducible $V$-modules
in this way.

For each $i=1,2,3$ let ${\cal S}_i$ 
be a finite set of inequivalent irreducible 
$V$-modules which is closed under the action of $G$
and let ${\cal L}_i$ be the direct sum of $V$-modules
in ${\cal S}_i$.
We have the decomposition ${\cal L}_i=
\bigoplus_{(j_i,\lambda_i)\in \Gamma_i}
W^i_{(j_i,\lambda_i)}\otimes M^i_{(j_i,\lambda_i)}$ as a
${\cal A}_{\alpha_i}(G,{\cal S}_i)\otimes V^G$-module.
Set 
${\cal I}=\bigoplus_{(L^1,L^2,L^3)\in {\cal S}_1\times{\cal S}_2\times{\cal S}_3}
I_V\binom{L^3}{L^1\ L^2}\otimes L^1\otimes L^2$, where 
$I_V\binom{L^3}{L^1\ L^2}$ is the set of 
all intertwining operators of type $\binom{L^3}{L^1\ L^2}$.
${\cal I}$ has a natural ${\cal A}_{\alpha_3}(G,{\cal S}_3)$-module structure.
For each $i=1,2,3$, fix $(j_i,\lambda_i)\in\Gamma_i$ and
nonzero 
$v^{10}\in M^1_{(j_1,\lambda_1)}, v^{20}\in M^1_{(j_2,\lambda_2)}$.
Set ${\cal I}_{(j_1,\lambda_1),(j_2,\lambda_2)}=
\mbox{Span}_{\C}\{f\otimes (w^1\otimes v^{10})\otimes 
(w^2\otimes v^{20})\in{\cal I}\ |\ 
w^1\in W^1_{(j_1,\lambda_1)}, 
w^2\in W^2_{(j_2,\lambda_2)}\}$.
${\cal I}_{(j_1,\lambda_1),(j_2,\lambda_2)}$ is 
a ${\cal A}_{\alpha_3}(G,{\cal S}_3)$-submodule of ${\cal I}$.
We will construct an injective linear map from 
the multiplicity space of $W^3_{(j_3,\lambda_3)}$ in 
${\cal I}_{(j_1,\lambda_1),(j_2,\lambda_2)}$ to the set of 
all intertwining operators for $V^G$ of  type $\binom{M^3_{(j_3,\lambda_3)}}
{M^1_{(j_1,\lambda_1)}\ M^2_{(j_2,\lambda_2)}}$.
Therefore, 
the fusion rule for $V^G$ of type $\binom{M^3_{(j_3,\lambda_3)}}
{M^1_{(j_1,\lambda_1)}\ M^2_{(j_2,\lambda_2)}}$ is greater than or equal to
the multiplicity of $W^3_{(j_3,\lambda_3)}$ in 
${\cal I}_{(j_1,\lambda_1),(j_2,\lambda_2)}$.

This paper is organized as follows.
In Sect.2 we first recall 
a construction of 
irreducible $V^G$-modules from irreducible $V$-modules in \cite{DY}.
We also recall the definitions of intertwining operators and fusion rules.
In Sect.3 we give a construction of 
intertwining operators for irreducible $V^G$-modules which occur as
submodules of irreducible $V$-modules.
In Sect.4 we apply the main result to determine some fusion rules for 
a vertex operator algebra ${\cal W}$ studied in \cite{DLTYY}.
In Appendix, we give some singular vectors in ${\cal W}$-modules used 
in Sect. 4.

\section{Preliminary}
We assume that the reader is familiar with the basic knowledge on
vertex operator algebras as presented in \cite{B}, \cite{FHL}, and \cite{FLM}. 

The following notation will be in force throughout the paper: $V=
(V, Y, \1,\omega)$ is a simple vertex operator algebra and $G$
a finite automorphism group of $V$.
For any $V$-module $L$, we always arrange the grading on
$L=\oplus_{n=0}^{\infty}L(n)$ so that $L(0)\neq 0$ if $L\neq 0$
by using  a grading shift.

\subsection{Irreducible $V^G$-modules constructed from irreducible $V$-modules}

In this subsection we review the results of Dong--Yamskulna\cite{DY}.
For a simple vertex operator algebra $V$, they showed that every irreducible $V$-module is 
a completely reducible $V^G$-module as a natural consequence
of a duality theorem of Schur-Weyl type.

Let $(L,Y_L)$ be an irreducible $V$-module and $a\in G$.
We define a new irreducible $V$-module $(L\circ a,Y_{L\circ a})$.
Here $L\circ a$ is equal to $L$ and $Y_{L\circ a}(u,z)=Y_L(au,z)$.
Note that $L\circ a$ is also an
irreducible $V$-module.
A set ${\cal S}$ of irreducible $V$-modules is called {\it stable}
if for any $L\in {\cal S}$ and $a\in G$ 
there exists $M\in{\cal S}$ such that $L\circ a\simeq M$. 

Now we take a finite $G$-stable set ${\cal S}$ consisting of 
inequivalent irreducible $V$-modules.
Let $L\in{\cal S}$ and $a\in G$. Then there exists $M\in{\cal S}$
such that $M\simeq L\circ a^{-1}$.
That is, there is a linear map $\phi(a,L) : L\rightarrow M$ satisfying the 
condition: $\phi(a,L)Y_L(v,z)=Y_M(av,z)\phi(a,L)$. By simplicity of $L$,
there exists $\alpha_L(b,a)\in\C^{\times}$ such that
$\phi(b,M)\phi(a,L)=\alpha_{L}(b,a)\phi(ba,L)$.
Moreover, for $a,b,c\in G$ we have
\[\alpha_L(c,ba)\alpha_L(b,a)=\alpha_M(c,b)\alpha_L(cb,a).\]
For $L\in{\cal S}$ and $a\in G$, we denote 
$M\in{\cal S}$ such that $L\circ a\simeq M$ by $L\cdot a$.

Define a vector space $\C{\cal S}=\bigoplus_{L\in {\cal S}}\C e(L)$ with a basis 
$e(L)$ for $L\in{\cal S}$. The space $\C{\cal S}$ is an associative algebra under
the product $e(L)e(M)=\delta_{L,M}e(L)$.
Let ${\cal U}(\C{\cal S})=\{\sum_{L\in{\cal S}}\lambda_Le(L)\ |\ 
\lambda_L\in\C^{\times}\}$ be the set of unit elements on $\C{\cal S}$.
${\cal U}(\C{\cal S})$ is a multiplicative right $G$-set by the action
$(\sum_{L\in{\cal S}}\lambda_{L}e(L))\cdot a=
\sum_{L\in{\cal S}}\lambda_{L}e(L\cdot a)$ for $a\in G$.
Set $\alpha(a,b)=\sum_{L\in{\cal S}}\alpha_{L}(a,b)e(L)$.
Then $(\alpha(a,b)\cdot c)\alpha(ab,c)=\alpha(a,bc)\alpha(b,c)$
hold for all $a,b,c\in G$.  So $\alpha : G\times G\rightarrow {\cal U}(\C{\cal S})$
is a $2$-cocycle.

Define the vector space
${\cal A}_{\alpha}(G,{\cal S})={\mathbb C}[G]
\otimes {\mathbb C}{\cal S}$ with a basis $a\otimes e(L)$ for
$a\in G$ and $L\in {\cal S}$ and a multiplication on it:
\[a\otimes e(L)\cdot b\otimes e(M)=
\alpha_{M}(a,b)ab\otimes e(L\cdot b)e(M).\]
Then ${\cal A}_{\alpha}(G,{\cal S})$ is an associative algebra
with the identity element 
$\sum_{L\in {\cal S}}1\otimes e(L)$.

We define an action of ${\cal A}_{\alpha}(G,{\cal S})$ on 
${\cal L}=\bigoplus_{L\in {\cal S}}L$ as follows:
For $L, M\in {\cal S},\ w\in M$ and $a\in G$ we set 
\begin{eqnarray*}
a\otimes e(L)\cdot w & = & \delta_{L,M}\phi(a,L)w.
\end{eqnarray*}
Note that the actions of ${\cal A}_{\alpha}(G,{\cal S})$
and $V^G$ on ${\cal L}$ commute with each other.

For each $L\in {\cal S}$ set $G_{L}
=\{a\in G\ |\ L\circ a\simeq L
\mbox{ as $V$-modules}\}$. 
Let ${\cal O}_{L}$
be the orbit of $L$ under the action of $G$ 
and let $G=\sqcup_{j=1}^{k}G_{L}g_j$ 
be a right coset decomposition with $g_1=1$.
Then ${\cal O}_{L}=\{L\cdot g_j\ |\ j=1,\ldots,k\}$ and 
$G_{L\cdot g_j}=g_j^{-1}G_{L}g_j$. We define
several subspaces of ${\cal A}_{\alpha}(G,{\cal S})$ by:
\[
\begin{array}{rcl}
{\displaystyle S(L)} & = & 
{\displaystyle 
\mbox{Span}_{\C}\{a\otimes e(L)\ |\ a\in G_{L}\}},\\
{\displaystyle D(L)} & = & 
{\displaystyle \mbox{Span}_{\C}\{a\otimes e(L)\ |\ a\in G\}} 
\mbox{  and} \\
{\displaystyle D({\cal O}_{L}}) & = & 
{\displaystyle \mbox{Span}_{\C}\{a\otimes 
e(L\cdot g_j)\ |\ j=1,\ldots,k, a\in G\}}.
\end{array}
\]
Decompose ${\cal S}$ into a disjoint union of orbits 
${\cal S}=\sqcup_{j\in J}{\cal O}_j$. Let $L^{(j)}$ be a 
representative element of ${\cal O}_j$. Then 
${\cal O}_j=\{L^{(j)}\cdot a\ |\ a\in G\}$
and ${\cal A}_{\alpha}(G, {\cal S})
=\bigoplus_{j\in J}D({\cal O}_{L^{(j)}})$.
We recall the following properties of 
${\cal A}_{\alpha}(G,{\cal S})$.

\begin{lemma}{\rm (\cite{DY}, Lemma 3.4)}
Let $L\in {\cal S}$ and $G=\sqcup_{j=1}^kG_{L}g_j$. Then
\begin{enumerate}
\item $S(L)$ is a subalgebra of ${\cal A}_{\alpha}(G,{\cal S})$
isomorphic to ${\mathbb C}^{\alpha_{L}}[G_{L}]$ where 
${\mathbb C}^{\alpha_{L}}[G_{L}]$ is the twisted 
group algebra with 
$2$-cocycle $\alpha_{L}$.
\item $D({\cal O}_{L})=\oplus_{j=1}^{k}D(L\cdot g_j)$ 
is a direct sum of left ideals.
\item Each $D({\cal O}_{L})$ is a two sided ideal of 
${\cal A}_{\alpha}(G, {\cal S})$ and 
${\cal A}_{\alpha}(G, {\cal S})
=\oplus_{j\in J}D({\cal O}_{L^{(j)}})$.
Moreover, $D({\cal O}_{L})$ has the identity element 
$\sum_{M\in {\cal O}_{L}}1\otimes e(M)$.
\end{enumerate}
\end{lemma}

\begin{lemma}{\rm (\cite {DY}, Theorem 3.6)}\label{Lemma : rep} 
\begin{enumerate}
\item $D({\cal O}_{L})$ is semisimple for all 
$L\in {\cal S}$
and the simple $D({\cal O}_{L})$-modules are precisely equal to 
$\mbox{\rm Ind}^{D(L)}_{S(L)}U=D(L)\otimes_{S(L)}U$ 
where $U$ ranges over the simple 
${\mathbb C}^{\alpha_{L}}[G_{L}]$-modules.
\item ${\cal A}_{\alpha}(G, {\cal S})$ 
is semisimple and simple ${\cal A}_{\alpha}(G, {\cal S})$-modules 
are precisely
$\mbox{\rm Ind}^{D(L^{(j)})}_{S(L^{(j)})}U$ 
where $U$ ranges over the simple 
${\mathbb C}^{\alpha_{L^{(j)}}}[G_{L^{(j)}}]$-modules 
and $j\in J$.
\end{enumerate}
\end{lemma}

For $L\in {\cal S}$ 
let $\Lambda_{G_{L}, \alpha_{L}}$ be the set of 
all irreducible characters $\lambda$ of 
${\mathbb C}^{\alpha_{L}}[G_{L}]$. We denote the 
corresponding
simple module by $U(L,{\lambda})$. Note that $L$ is a
semisimple ${\mathbb C}^{\alpha_{L}}[G_{L}]$-module.
Let $L^{\lambda}$ be the sum of simple 
${\mathbb C}^{\alpha_{L}}[G_{L}]$-submodules of $L$ 
isomorphic to $U(L,{\lambda})$. Then 
$L=\oplus_{\lambda\in \Lambda_{G_{L},\alpha_{L}}}L^{\lambda}$.
Moreover $L^{\lambda}=U(L,{\lambda})\otimes L_{\lambda}$
where $L_{\lambda}=
\mbox{Hom}_{{\mathbb C}^{\alpha_{L}}[G_{L}]}
(U(L,{\lambda}),L)$ is the multiplicity space of $U(L,{\lambda})$ 
in $L$.
We can realize $L_{\lambda}$ as a subspace of $L$
in the following way: Let $w\in U(L,{\lambda})$ be a 
fixed nonzero
vector. Then we can identify 
$\mbox{Hom}_{{\mathbb C}^{\alpha_{L}}[G_{L}]}
(U(L,{\lambda}),L)$  with the subspace 
\[\{f(w)\ |\ f\in 
\mbox{Hom}_{{\mathbb C}^{\alpha_{L}}[G_{L}]}
(U(L,{\lambda}),L)\}\]
of $L^{\lambda}$. Note that the actions of 
${\mathbb C}^{\alpha_{L}}[G_{L}]$
and $V^{G_{L}}$ on $L$ commute with each other. So
$L^{\lambda}$ and $L_{\lambda}$ are ordinary 
$V^{G_{L}}$-modules. 
Furthermore, 
$L^{\lambda}$ and $L_{\lambda}$ are ordinary $V^{G}$-modules. 

For convenience, we set 
$G_j=G_{L^{(j)}}, 
\Lambda_j=\Lambda_{L^{(j)},\alpha_{L^{(j)}}}$ and
$U_{(j, \lambda)}=U(L^{(j)},{\lambda})$ for $j\in J$ and
$\lambda\in \Lambda_j$.
We denote by $\Gamma$ the set 
$\{(j,\lambda)\ |\ j\in J, \lambda\in\Lambda_j\}$.
We have a decomposition
\[L^{(j)}=\bigoplus_{\lambda\in \Lambda_j}
U_{(j,\lambda)}\otimes M_{(j,\lambda)}\]
as a ${\mathbb C}^{\alpha_{L^{(j)}}}[G_{j}]
\otimes V^{G_{j}}$-module.
We also have
\[{\cal L}=\bigoplus_{(j,\lambda)\in \Gamma}
\mbox{Ind}^{D({L^{(j)}})}_{S({L^{(j)}})}
U_{(j,\lambda)}\otimes M_{(j,\lambda)}\]
as a ${\cal A}_{\alpha}(G, {\cal S})\otimes V^G$-module.
 
For $(j,\lambda)\in \Gamma$ we set $W_{(j,\lambda)}=
\mbox{Ind}^{D(L^{(j)})}_{S(L^{(j)})}U_{(j, \lambda)}$.
Then $W_{(j,\lambda)}$ forms a complete list of simple 
${\cal A}_{\alpha}(G, {\cal S})$-modules by
Lemma \ref{Lemma : rep}.

A duality theorem of Schur--Weyl type holds.

\begin{theorem}\label{theorem : DY}{\rm (\cite{DY} Theorem 6.14)}
As a ${\cal A}_{\alpha}(G,{\cal S})\otimes V^G$-module,
\[{\cal L}=\bigoplus_{(j,\lambda)\in \Gamma}
W_{(j,\lambda)}\otimes M_{(j,\lambda)}.\]
Moreover,
\begin{enumerate}
\item Each $M_{(j,\lambda)}$ is a nonzero irreducible $V^G$-module.
\item $M_{(j_1,\lambda_1)}$ and $M_{(j_2,\lambda_2)}$
are isomorphic $V^G$-modules if and only if 
$(j_1,\lambda_1)=(j_2,\lambda_2)$.
\end{enumerate}
\end{theorem}

\subsection{Intertwining operators and fusion rules}

We recall the definition of intertwining operators for $V$-modules which is 
introduced in \cite{FHL}.

\begin{definition}{\rm
Let $L^i\ (i=1,2,3)$ be $V$-modules.
An intertwining operator of type 
$\binom{L^3}{L^1\ L^2}$ is a linear map 
\[\begin{array}{rrcl}
I(\cdot,z) : & L^1 & \rightarrow & \mbox{Hom}_{\C}(L^2,L^3)\{z\},\\
& v & \mapsto & I(v,z)=\sum_{\gamma\in \C}v_{\gamma}^{-\gamma-1},
\end{array}\]
which satisfies the following conditions:
Let $u\in V, v\in L^1$, and $w\in L^2$.
\begin{enumerate}
\item For any fixed 
$\gamma\in\C$, $v_{\gamma+n}w=0$ for $n\in\Z$ sufficiently large.
\item $\displaystyle{I(L(-1)v,z)=\frac{d}{dz}I(v,z)}$.
\item 
\begin{eqnarray*}
& & z_0^{-1}\delta(\frac{z_1-z_2}{z_0})Y_{L^3}(u,z_1)I(v,z_2)w-
z_0^{-1}\delta(\frac{-z_2+z_1}{z_0})I(v,z_2)Y_{L^2}(u,z_1)w\\
& = & \ 
z_2^{-1}\delta(\frac{z_1-z_0}{z_2})I(Y_{L^1}(u,z_0)v,z_2)w.
\end{eqnarray*}
\end{enumerate}
We denote by $I_V\binom{L^3}{L^1\ L^2}$ the set of all 
intertwining operators of type $\binom{L^3}{L^1\ L^2}$.
The dimension of $I_V\binom{L^3}{L^1\ L^2}$
is called the {\it fusion rule} of type $\binom{L^3}{L^1\ L^2}$.
}
\end{definition}

For a $V$-module $L=\oplus_{n=0}^{\infty}L(n)$, it is shown in \cite[Theorem 5.2.1.]{FHL} that 
the graded vector space $L^{\prime}=\oplus_{n=0}^{\infty}L(n)^{*}$ carries the structure of a $V$-module,
where $L(n)^{*}=\mbox{Hom}_{\C}(L(n),\C)$.
$L^{\prime}$ is called the {\it contragredient module} of $L$.

The fusion rules have some symmetries.
\begin{lemma}\label{lemma : sym}{\rm (\cite{FHL}, Proposition 5.4.7 and 5.5.2)}
Let $L^i\ (i=1,2,3)$ be $V$-modules.
Then
\[
\dim_{\C}I_V\binom{L^3}{L^1\ L^2}=
\dim_{\C}I_V\binom{L^3}{L^2\ L^1}=
\dim_{\C}I_V\binom{(L^2)^{\prime}}{L^1\ (L^3)^{\prime}}. 
\]
\end{lemma}
Let $L^1$ and $L^2$ be irreducible $V$-modules.
We use a notation 
$L^1\times L^2=\sum_{L^3}\dim_{\C}I_{V}\binom{L^3}{L^1\ L^2}L^3$
to represent the fusion rules, where $L^3$ ranges over the 
irreducible $V$-modules. Note that 
$L^1\times L^2=L^2\times L^1$ by Lemma \ref{lemma : sym}.

\section{Intertwining operators for irreducible $V^G$-modules
which occur as submodules of irreducible $V$-modules}
In this section we give a construction of intertwining operators
for irreducible $V^G$-modules which occur as submodules of irreducible $V$-modules
by using  intertwining operators for $V$.

Let ${\cal S}_i\ (i=1,2,3)$ be 
finite $G$-stable sets consisting of 
inequivalent irreducible $V$-modules.
Set ${\cal L}_{i}=\bigoplus_{L\in{\cal S}_i}L$.
For $L^i\in {\cal S}_i$ and $a\in G$,
$\phi_i(a,L^i) : L^i\rightarrow L^i\cdot a^{-1}$ 
denote the fixed $V$-module isomorphisms.
For $L^i\in{\cal S}_i$ and $a,b\in G$, $\alpha_{i,L^i}\in\C^{\times}$
denote nonzero complex numbers such that 
$\phi_i(b,L^i\cdot a^{-1})\phi_i(a,L^i)=\alpha_{i,L^i}(b,a)\phi(ba,L^i)$.
Set $\alpha_i(a,b)=\sum_{L^i\in{\cal S}_i}\alpha_{i,L^i}(a,b)e(L^i)$.
Let ${\cal S}_i=\sqcup_{j\in J_i}{\cal O}_j$ be the orbit decompositions
under the action of $G$ and set $\Gamma_i=\{(j_i,\lambda_i)\ |\ 
j_i\in J_i, \lambda_i\in \Lambda_{j_i}\}$.

For $f\in I_V\binom{L^3}{L^1\ L^2}$ and $a\in G$, we define 
${}_{a}f\in I_V\binom{L^3\cdot a^{-1}}
{L^1\cdot a^{-1}\ L^2\cdot a^{-1}}$ as follows:
For $v\in L^1$ we set 
\[{}_af(v,z)=\phi_3(a,L^3)
f(\phi_1(a,L^1)^{-1}v, z)\phi_2(a,L^2)^{-1}.\]
Set
\[{\cal I}=\bigoplus_{(L^1,L^2,L^3)\in
{\cal S}_1\times{\cal S}_2\times{\cal S}_3}
I_V\binom{L^3}{L^1\ L^2}\otimes_{\C}L^1\otimes_{\C}L^2. \]
We define an action of ${\cal A}_{\alpha_3}(G,{\cal S})$ on ${\cal I}$ as follows:
Let $L^i\in {\cal S}_i\ (i=1,2,3)$. 
For $a\otimes e(M)\in{\cal A}_{\alpha_3}(G,{\cal S}),
v\in L^1, w\in L^2$, and $f\in 
I_V\binom{L^3}{L^1\ L^2}$, we set
\begin{eqnarray*}
(a\otimes e(M))\cdot(f\otimes v\otimes w) & = &
\delta_{M,L^3}\cdot{}_{a}f\otimes \phi_1(a,L^1)v\otimes \phi_2(a,L^2)w\\
& & \in I_V\binom{L^3\cdot a^{-1}}{L^1\cdot a^{-1}\ L^2\cdot a^{-1}}\otimes 
L^1\cdot a^{-1}\otimes L^2\cdot a^{-1}.\end{eqnarray*}

We define a map $\Psi : {\cal I}\rightarrow {\cal L}_3\{z\}$ by
\[\Psi(f\otimes v^{1}\otimes v^2)=f(v^1,z)v^2\]
for $v^1\in L^1, v^2\in L^2$, and $f\in I_V\binom{L^3}{L^1\ L^2}$
,where $L^i\in{\cal S}_i\ (i=1,2,3)$.
Note that $\Psi$ is a ${\cal A}_{\alpha_3}(G,{\cal S}_3)$-module homomorphism.

\begin{lemma}\label{lemma : indep}The map 
$\Psi :
{\cal I}\rightarrow {\cal L}_3\{z\}$ is injective.
\end{lemma}
\begin{proof}
We use the same method that was used in the proof of Lemma 3.1. of \cite{DM}.
Assume false. Then there is a nonzero $X\in \mbox{Ker}\Psi$.
Since ${\cal L}_3=\oplus_{L\in{\cal S}_3}L$,
we may assume $X=\sum_{i,j}f^{ij}\otimes v^{1i}\otimes v^{2j}$, where
$v^{1i}\in L^{1i} (i=1,\ldots,l_1)$ are linearly independent homogeneous vectors in ${\cal L}_1$,
$v^{2j}\in L^{2j} (j=1,\ldots,l_2)$ are linearly independent homogeneous vectors in ${\cal L}_2$,
$f^{ij}\in I_V\binom{L^{3}}{L^{1i}\ L^{2j}}$, $L^{1i}\in {\cal S}_1$,
$L^{2j}\in {\cal S}_2$, and $L^3\in{\cal S}_3$.
We may also assume $f^{11}\otimes v^{11}\otimes v^{21}$ is nonzero.
Since $\sum_{i,j}f^{ij}(v^{1i},z)v^{2j}=0$,
for all $u\in V$ we have
\[\sum_{i,j}Y_{{\cal L}_3}(u,z_1)f^{ij}(v^{1i},z)v^{2j}=0.\]
Using the associativity of intertwining operators \cite[Proposition 11.5]{DL},
\begin{eqnarray}\label{formula : ass}
\sum_{i,j}f^{ij}(Y_{{\cal L}_1}(u,z_1)v^{1i},z)v^{2j} & = & 0.
\end{eqnarray}

We denote $Y_{{\cal L}_1}(u,z_1)=\sum_{n\in\Z}u^{{\cal L}_1}_{n}z_1^{-n-1}$.
Fix $N\in\Z$ such that 
$v^{1i}\in \oplus_{n=0}^{N}L^{1i}(n)$ for all $i$.
Since ${\cal S}_1$ consists of inequivalent irreducible $V$-modules,
the linear map $\sigma_N : V\rightarrow \oplus_{L\in S_1}\oplus_{m=0}^{N}\mbox{End}_{\C}L(m)$ 
defined by 
$\sigma_{N}(u)=u^{{\cal L}_1}_{\mmbox{wt}u-1}$ for
homogeneous $u\in V$ is an epimorphism by Lemma 6.13 of \cite{DY}.
So there exists $u^1\in V$ such that $\sigma_N(u^1)v^{1i}=\delta_{1,i}v^{11}$.
{From} formula (\ref{formula : ass}), we have
\begin{eqnarray}\label{formula : result}
0 & =& \sum_{i,j}f^{ij}(\sigma_N(u^1)v^{1i},z)v^{2j}
=\sum_{j}f^{1j}(v^{11},z)v^{2j}.
\end{eqnarray}
Therefore, for all $u\in V$ we have
\[\sum_{j}Y_{{\cal L}_3}(u,z_1)
f^{1j}(v^{11},z)v^{2j}=0.\]
Using the commutativity of intertwining operators \cite[Proposition 11.4]{DL},
\[\sum_{j}f^{1j}(v^{11},z)Y_{{\cal L}_2}(u,z_1)v^{2j}=0.\]
Since ${\cal S}_2$ consists of inequivalent irreducible $V$-modules, we have
\begin{eqnarray*}
f^{11}(v^{11},z)v^{21} & = & 0
\end{eqnarray*}
for the same reason to obtain formula (\ref{formula : result}).
So $f^{11}\otimes v^{11}\otimes v^{21}\in \mbox{Ker}\Psi$.
Since $L^{11}$ and $L^{21}$ are irreducible $V$-modules and 
$f^{11}\otimes v^{11}\otimes v^{21}$ is nonzero, 
this contradicts Proposition 11.9 of \cite{DL}. 
\end{proof}
We have the decomposition of each ${\cal L}_i$ as 
a ${\cal A}_{\alpha_i}(G,{\cal S}_i)$-module in Theorem \ref{theorem : DY}:
\begin{eqnarray*}
{\cal L}_i & = & \bigoplus_{(j_i,\lambda_i)\in \Gamma_i}
W^i_{(j_i,\lambda_i)}\otimes M^i_{(j_i,\lambda_i)}.
\end{eqnarray*}
For $i=1,2$ let $(j_i,\lambda_i)\in \Gamma_i$ and let $v^{i0}\in M^i_{(j_i,\lambda_i)}$.
Set
\begin{eqnarray*}
\lefteqn{{\cal I}_{(j_1,\lambda_1),(j_2,\lambda_2)}(v^{10},v^{20})}
\\
& = & \mbox{Span}_{\C}\{f\otimes (w^1\otimes v^{10})\otimes (w^2\otimes v^{20})\in{\cal I}\ |\ 
w^1\in W^1_{(j_1,\lambda_1)}, 
w^2\in W^2_{(j_2,\lambda_2)}\}.
\end{eqnarray*}
${\cal I}_{(j_1,\lambda_1),(j_2,\lambda_2)}(v^{10},v^{20})$
is a ${\cal A}_{\alpha_3}(G,{\cal S}_3)$-submodule of ${\cal I}$.
For any nonzero $v^{10},v^1\in M^1_{(j_1,\lambda_1)}$ and nonzero
$v^{20},v^{2}\in M^2_{(j_2,\lambda_2)}$,
${\cal I}_{(j_1,\lambda_1),(j_2,\lambda_2)}(v^{10},v^{20})$
and ${\cal I}_{(j_1,\lambda_1),(j_2,\lambda_2)}(v^{1},v^{2})$
are isomorphic ${\cal A}_{\alpha_3}(G,{\cal S}_3)$-modules.

\begin{theorem}\label{theorem : main}{\rm
Fix a nonzero $v^{10}\in M^1_{(j_1,\lambda_1)}$ and a nonzero $v^{20}\in M^2_{(j_2,\lambda_2)}$.

For any
$((j_1,\lambda_1),(j_2,\lambda_2),(j_3,\lambda_3))\in 
\Gamma_1\times\Gamma_2\times\Gamma_3$, there exists an
injective linear map from 
$\mbox{Hom}_{{\cal A}_{\alpha_3}(G,{\cal S}_3)}(W^3_{(j_3,\lambda_3)},
{\cal I}_{(j_1,\lambda_1),(j_2,\lambda_2)}(v^{10},v^{20}))$
to $I_{V^G}\binom{M^3_{(j_3,\lambda_3)}}
{M^1_{(j_1,\lambda_1)}\ M^2_{(j_2,\lambda_2)}}$.
In particular,
\[\dim_{\mathbb C}I_{V^G}\binom{M^3_{(j_3,\lambda_3)}}
{M^1_{(j_1,\lambda_1)}\ M^2_{(j_2,\lambda_2)}}\geq 
\dim_{\mathbb C}\mbox{Hom}_{{\cal A}_{\alpha_3}(G,{\cal S}_3)}(W^3_{(j_3,\lambda_3)},
{\cal I}_{(j_1,\lambda_1),(j_2,\lambda_2)}(v^{10},v^{20})).\]
}
\end{theorem}
\begin{proof}
For convenience, we set 
$W^i=W^i_{(j_i,\lambda_i)}, M^i=M^i_{(j_i,\lambda_i)}\ (i=1,2,3),
{\cal A}_3={\cal A}_{\alpha_3}(G,{\cal S}_3)$, and 
${\cal I}_{12}(v^{10},v^{20})={\cal I}_{(j_1,\lambda_1),(j_2,\lambda_2)}(v^{10},v^{20})$.
Fix a nonzero $w^{30}\in W^3$.
Let $F\in \mbox{Hom}_{{\cal A}_3}(W^3,
{\cal I}_{12}(v^{10},v^{20}))$. 
We shall define $\Phi(F)\in I_{V^G}\binom{M^3}
{M^1\ M^2}$.

For $v^1\in M^1, v^2\in M^2$,
we set \[H(v^1,v^2)\in 
\mbox{Hom}_{{\cal A}_3}(W^3,
{\cal I}_{12}(v^1,v^2))\]
as follows:
For $w^3\in W^3$, 
let $F(w^3)=\sum_{i}f^i\otimes (w^{1i}\otimes v^{10})
\otimes (w^{2i}\otimes v^{20})\in 
{\cal I}_{12}(v^{10},v^{20})$, where
$w^{1i}\otimes v^{10}\in L^{1i}, w^{2i}\otimes v^{20}\in L^{2i}$,
$f^i\in I_V\binom{L^{3i}}{L^{1i}\ L^{2i}}$ and 
$L^{ji}\in{\cal S}_j\ (j=1,2,3)$.
Note that $w^{ji}\otimes v^{j}\in L^{ji}$
by the definition of $M^j$ for $j=1,2$.
We define 
\begin{eqnarray*}
H(v^1,v^2)(w^3) & = & \sum_{i}f^i\otimes (w^{1i}\otimes v^1)\otimes (w^{2i}\otimes v^2)
\end{eqnarray*}
It is clear that $H(v^1,v^2)\in\mbox{Hom}_{{\cal A}_3}(W^3,{\cal I}_{12}(v^1,v^2))$. 
Since the map
$\Psi : {\cal I}\rightarrow {\cal L}_3\{z\}$ 
is a ${\cal A}_3$-module homomorphism, 
the map $w^3\mapsto \Psi(H(v^1,v^2)(w^3))$
is a ${\cal A}_3$-module homomorphism from $W^3$ to ${\cal L}_3\{z\}$.
%
So $\Psi(H(v^1,v^2)(W^3))$ is a ${\cal A}_3$-submodule of ($W^3\otimes M^3)\{z\}$.
Let $w^{3,1},w^{3,2},\ldots,w^{3,\dim_{\C}W^3}$
be a basis of $W^3$. Since $W^3$ is an irreducible ${\cal A}_3$-module,
there exists $a\in {\cal A}_3$ such that $aw^{3,i}=\delta_{1,i}w^{3,1}$. 
Let $\Psi(H(v^1,v^2)(w^{3,1}))=\sum_{i}w^{3,i}\otimes p^i$, where 
$p_i\in M^3\{z\}$. Then
\begin{eqnarray*}
\lefteqn{\Psi(H(v^1,v^2)(w^{3,1})) = \Psi(H(v^1,v^2)(aw^{3,1}))=a\Psi(H(v^1,v^2)(w^{3,1}))}\\
& = & a\sum_{i}w^{3,i}\otimes p^i=\sum_{i}(aw^{3,i})\otimes p^i=w^{3,1}\otimes p^1.
\end{eqnarray*}
So $\Psi(H(v^1,v^2)(w^3))\in (w^3\otimes M^3)\{z\}$ for all $w^3\in W^3$.
We hence have an unique $\Phi(F)(v^1,z)v^2\in M^3\{z\}$ such that 
\begin{eqnarray}
w^{30}\otimes \Phi(F)(v^1,z)v^2 
& = & \Psi(H(v^1,v^2)(w^{30})).\label{eqn : inter}
\end{eqnarray}
Since $f^i$ are intertwining operators,
we have $\Phi(F)\in I_{V^G}\binom{M^3}
{M^1\ M^2}$
from formula (\ref{eqn : inter}).
 
We will show that $\Phi$ is injective.
Suppose $\Phi(F)=0$. Then 
\[
0 = w^{30}\otimes \Phi(F)(v^{10},z)v^{20}=\Psi(F(w^{30})).
\]
Since $\Psi$ is injective by Lemma \ref{lemma : indep},
$F(w^{30})=0$.
Since $W^3$ is an irreducible
${\cal A}_3$-module and 
$F\in \mbox{Hom}_{{\cal A}_3}(W^3, {\cal I}_{12}(v^{10},v^{20}))$, $F=0$.
\end{proof}
Let $\C G$ be the group algebra of $G$ and $\mbox{Irr}G$
the set of all irreducible characters of $G$.
We set ${\cal S}_i=\{V\}\ (i=1,2,3)$ in Theorem \ref{theorem : main}.
Then ${\cal A}_{\alpha_i}(G,{\cal S}_i)=\C G$ 
and $\Gamma_i=\mbox{Irr}G$.
Note that $\dim_{\C}I_V\binom{V}{V\ V}=1$ since $V$ is simple.
In this case we have the following result:
\begin{corollary}\label{corollary : group}
Let $\chi_i\in\mbox{\rm Irr}G\ (i=1,2,3)$. Then
\[\dim_{\mathbb C}I_{V^G}\binom{V_{\chi_3}}
{V_{\chi_1}\ V_{\chi_2}}\geq 
\dim_{\C}\mbox{\rm Hom}_{\C G}(W_{\chi_3},
W_{\chi_1}\otimes_{\C}W_{\chi_2}).\]
\end{corollary}
In \cite[Section 3]{DVVV}, it is conjectured that if $V$ is rational then
for all $\chi_1,\chi_2\in \mbox{\rm Irr}G$,
\begin{eqnarray*}
V_{\chi_1}\times  V_{\chi_2} & = & 
\sum_{\chi_3\in\mbox{\rm \scriptsize Irr}G}\dim_{\C}\mbox{\rm Hom}_{\C G}(W_{\chi_3},
W_{\chi_1}\otimes_{\C}W_{\chi_2})V_{\chi_3}.
\end{eqnarray*}
The conjecture implies that if $V$ is rational then
the representation algebra of the finite group $G$ is always 
realized as a subalgebra of the fusion algebra of $V^G$.

\section{An application}
In \cite{DLTYY} we studied a vertex operator algebra
${\cal W}$ which is a realization of an algebra denoted by $[Z^{(5)}_3]$ in \cite{FaZ}.
${\cal W}$ is a fixed point subalgebra of a vertex operator algebra $M_k^0$.
It is expected that the $\Z_3$ symmetry of ${\cal W}$ affords $3B$ elements 
of the Monster simple group \cite{KLY2}.
For $M_k^0$, the irreducible modules are classified and 
the fusion rules are determined in \cite{LY}.
Let $L^i\ (i=1,2,3)$ be irreducible ${\cal W}$-modules such 
that $L^1$ and $L^2$ occur as submodules of irreducible $M^0_k$-modules.
In this section we determine the 
fusion rule of type $\binom{L^3}{L^1\ L^2}$
by using  Theorem \ref{theorem : main}. 

\subsection{Subalgebra $M_k^0$ of $V_{\sqrt{2}A_2}$}
In this subsection we review some properties of $M_k^0$ in \cite{LY}.
Let $A_2$ be the ordinary root lattice of type $A_2$
and $V_{\sqrt{2}A_2}$ the lattice vertex operator algebra
associated with  $\sqrt{2}A_2$.
Let $\alpha_1,\alpha_2$ be the simple roots of type $A_2$ and set 
$\alpha_0=-\alpha_1-\alpha_2$. 

For basic definitions concerning lattice vertex operator algebras we refer 
to \cite{DL} and \cite{FLM}.
Our notation for the lattice vertex operator algebra is standard 
\cite{FLM}.
In particular, $\h=\C\otimes_{\Z}\sqrt{2}A_2$ is an abelian Lie algebra,
$\hat{\h}=\h\otimes \C[t,t^{-1}]\oplus \C c$ is the corresponding affine Lie algebra,
$M(1)=\C[\alpha(n)\ ;\ \alpha\in\h, n<0\}$, where 
$\alpha(n)=\alpha\otimes t^n$, is the unique irreducible 
$\hat{\h}$-module such that $\alpha(n)1=0$ for all $\alpha\in\h, n>0$, and $c=1$.
As a vector space $V_{\sqrt{2}A_2}=M(1)\otimes \C[\sqrt{2}A_2]$ and 
for each $v\in V_{\sqrt{2}A_2}$, a vertex operator $Y(v,z)=\sum_{n\in\Z}v_nz^{-n-1}
\in\mbox{End}(V_{\sqrt{2}A_2})[[z,z^{-1}]]$ is defined.
The vector ${\bf 1}=1\otimes 1$ is called the vacuum vector.
We use the symbol $e^{\alpha}, \alpha\in \sqrt{2}A_2$ to denote a basis of 
$\C[\sqrt{2}A_2]$.

There exists an isometry $\tau$ of $\sqrt{2}A_2$
such that $\tau(\sqrt{2}\alpha_1)=\sqrt{2}\alpha_2$ and
$\tau(\sqrt{2}\alpha_2)=\sqrt{2}\alpha_0$.
The isometory $\tau$ lifts naturally to an automorphism of $V_{\sqrt{2}A_2}$:
\begin{eqnarray*}
\alpha^1(-n_1)\cdots\alpha^k(-n_k)e^{\beta}
& \mapsto & 
(\tau\alpha^1)(-n_1)\cdots(\tau\alpha^k)(-n_k)e^{\tau\beta}.
\end{eqnarray*}
By abuse of notation, we denote it by $\tau$.
Let $G$ be the cyclic group generated by $\tau$.
Set 
\begin{eqnarray*}
\omega^3 & = &
\frac{1}{15}(\alpha_1(-1)^2+\alpha_2(-1)^2+\alpha_0(-1)^2)\\
& & {}+
\frac{1}{10}
(e^{\sqrt{2}\alpha_1}+e^{-\sqrt{2}\alpha_1}+
e^{\sqrt{2}\alpha_2}+e^{-\sqrt{2}\alpha_2}+
e^{\sqrt{2}\alpha_0}+e^{-\sqrt{2}\alpha_0})
\end{eqnarray*}
and $M_k^0=\{v\in V_{\sqrt{2}A_2}\ |\ (\omega^3)_1v=0\}$.
Since $\tau\omega^3=\omega^3$, $M_k^0$ is invariant under the action of $\tau$.
$\tau$ is an automorphism group of $M_k^0$ of order $3$ by \cite[Theorem 2.1]{DLTYY}. 
Set
\[L_0=L,\quad L_a=\frac{\sqrt{2}\alpha_2}{2}+L,\quad
L_b=\frac{\sqrt{2}\alpha_0}{2}+L,\quad L_c=\frac{\sqrt{2}\alpha_1}{2}+L
\]
and
\begin{eqnarray*}
M_{k}^{i} & = & \{ v\in V_{L_i}\,|\,(\omega^3)_1 v=0\},\\
W_{k}^{i} & = & \{ v\in V_{L_i}\,|\,(\omega^3)_1 v=
\frac{2}{5}v\} ,\qquad \text{ for } i=0,a,b,c.
\end{eqnarray*}
It is shown in \cite{LY} 
that $\{M^i_k,W^i_k\ |\ i=0,a,b,c\}$ is the set of all irreducible $M_k^0$-modules
and the fusion rules are determined.

\subsection{Subalgebra ${\cal W}$ in $M^{\tau}$}
We denote by ${\cal W}$ the subalgebra $(M_k^0)^{\tau}$ of fixed points of $\tau$ in 
$M_k^0$.
We recall some properties of ${\cal W}$ in \cite{DLTYY}.
${\cal W}$ is generated by the Virasoro element $\omega$ and 
an element $J$ of weight $3$. 
Let 
$Y(\omega,z)=\sum_{n\in \Z}L(n)z^{-n-2}$ and
$Y(J,z)=\sum_{n\in \Z}J(n)z^{-n-3}$.
They satisfy the following commutation relations:
\begin{eqnarray}\label{eqn : com}
[L(m),\,L(n)] & = & (m-n)L(m+n)+\frac{m^3-m}{12}\cdot\frac{6}{5}
\cdot\delta_{m+n,0},\nonumber\\
{[}L(m),\,J(n)] & = & (2m-n)J(m+n),\nonumber\\
{[}J(m),\,J(n)] & = & (m-n)(22(m+n+2)(m+n+3) + 35(m+2)(n+2))L(m+n)\nonumber\\
& & {}-120(m-n)\Big( \sum_{k \le -2} L(k)L(m+n-k) +
\sum_{k \ge -1}L(m+n-k)L(k) \Big)\nonumber\\
& & {}-\frac{7}{10}m(m^2-1)(m^2-4)\delta_{m+n,0}.
\end{eqnarray}
${\cal W}$ has exactly $20$ irreducible modules.
$8$ irreducible ${\cal W}$-modules occur as submodules of irreducible $M_k^0$-modules.
We introduce those $8$ ${\cal W}$-modules.
$\tau$ acts on the irreducible $M_k^0$-modules as follows:
\[
\begin{array}{l}
M^0_k\circ\tau\simeq M^0_k, W^0_k\circ\tau\simeq W^0_k,\\
M^a_k\circ\tau\simeq M^c_k, M^c_k\circ\tau\simeq M^b_k, M^b_k\circ\tau\simeq M^a_k,\\
W^a_k\circ\tau\simeq W^c_k, W^c_k\circ\tau\simeq W^b_k,\mbox{ and }
W^b_k\circ\tau\simeq W^a_k.
\end{array}
\]
So $\{M^0_k\}, \{W^0_k\}, \{M^a_k,M^b_k,M^c_k\}$,
and $\{W^a_k,W^b_k,W^c_k\}$ are $G$-stable sets.
The automorphism $\tau$ of $V_{\sqrt{2}A_2}$ fixes $\omega^3$ and
so $W_k^0$ is invariant under $\tau$. Hence we can take $\tau$ as 
$\phi(\tau,W_k^0)$ in Section 2.1.
For these $G$-stable sets, we can take the $2$-cocycles $\alpha$ in Section 2.1
to be trivial.
Let $\xi=e^{2\pi\sqrt{-1}/3}$.
We set 
$M_k^{0 (i)}=\{v\in M^0_k\ |\ \tau v=\xi^i v\}$ and 
$W_k^{0 (i)}=\{v\in W^0_k\ |\ \tau v=\xi^i v\}$ for $i\in\Z$.
Note that ${\cal W}=M_k^{0 (0)}$.
By Theorem \ref{theorem : DY}, we have
$M_k^{0 (i)},W_k^{0 (i)},\ (i=0,1,2)\ M^a_k$, and $W^a_k$ are inequivalent 
irreducible ${\cal W}$-modules. Moreover,
$M^a_k\simeq M^b_k\simeq M^c_k$ and $W^a_k\simeq W^b_k\simeq W^c_k$ 
as ${\cal W}$-modules. The contragredient modules of these 
${\cal W}$-modules are
\[\begin{array}{ll}
(M_k^{0 (i)})^{\prime}\simeq M_k^{0 (-i)}, \  
(W_k^{0 (i)})^{\prime}\simeq W_k^{0 (-i)}, & (i=0,1,2),\\
(M_k^a)^{\prime} \simeq M_k^a \mbox{ and } (W_k^a)^{\prime} \simeq W_k^a.
\end{array}
\]

All the other irreducible ${\cal W}$-modules 
occur as submodules of irreducible $\tau^i$-twisted $M^0_k$-modules for $i=1,2$.

\subsection{An upper bound for the fusion rule}

We review some notations and formulas for the Zhu algebra $A(V)$
of an arbitrary vertex operator algebra $V$ in \cite{Z}
and the $A(V)$-bimodule $A(L)$ of an arbitrary $V$-module $L$ in \cite{FZ} and \cite{L}.
For $u,v\in V$ with $u$ being homogeneous, define two bilinear operations
\begin{eqnarray*}
u*v & = & \mbox{Res}_{z}(Y(u,z)v\frac{(1+z)^{\mmbox{wt}u}}{z}),\\
u\circ v & = & \mbox{Res}_{z}(Y(u,z)v\frac{(1+z)^{\mmbox{wt}u}}{z^2}).
\end{eqnarray*}
We extend $*$ and $\circ$ for arbitrary $u, v\in V$ by linearity. 
Let $O(V)$ be the subspace of $V$ spanned by all $u\circ v$ for $u\in V, v\in L$.
Set $A(V)=V/O(V)$.
By \cite[Theorem 2.1.1.]{Z}, $O(V)$ is a two-sided ideal with respect to the operation
$*$ and $(A(V),*)$ is an associative algebra with identity $1+O(V)$.
For every $V$-module $N$, $N(0)$ is a left $A(V)$-module.

Let $L$ be a $V$-module.
For $u\in V, v\in L$ with $u$ being homogeneous, define three bilinear operations
\begin{eqnarray*}
u*v & = & \mbox{Res}_{z}(Y(u,z)v\frac{(1+z)^{\mmbox{wt}u}}{z}),\\
v*u & = & \mbox{Res}_{z}(Y(u,z)v\frac{(1+z)^{\mmbox{wt}u-1}}{z}),\\
u\circ v & = & \mbox{Res}_{z}(Y(u,z)v\frac{(1+z)^{\mmbox{wt}u}}{z^2}).
\end{eqnarray*}
We extend $*$ and $\circ$ for arbitrary $u\in V, v\in L$ by linearity. 
Let $O(L)$ be the subspace of $L$ spanned by all $u\circ v$ for $u\in V, v\in L$.
By \cite[Theorem 1.5.1.]{FZ}, $O(L)$ is a two-sided ideal with respect to the operation
$*$. Thus it induces an operation on $A(L)=L/O(L)$.
Denote by $[v]$ the image of $v\in L$ in $A(L)$.
$A(L)$ is a $A(V)$-bimodule under the operation $*$.
Using $A(L)$, we have an upper bound for every fusion rule.

\begin{lemma}\label{lemma : ub}{\rm (\cite{L} Proposition 2.10.)}
Let $L^i=\bigoplus_{n=0}^{\infty}L^{i}(n)\ (i=1,2,3)$ be 
irreducible $V$-modules. Then
\[\dim_{\C}I_V\binom{L^3}{L^1\ L^2}\leq \dim_{\C}\mbox{\rm Hom}_{\C}
\Big(L^3(0)^{*}
\otimes_{A(V)}A(L^1)\otimes_{A(V)}L^2(0)\Big).\]
\end{lemma}

\subsection{The fusion rules for irreducible ${\cal W}$-modules 
which occur as submodules of irreducible $M_k^0$-modules}

Let $L^i\ (i=1,2,3)$ be irreducible ${\cal W}$-modules such that
$L^1$ and $L^2$ occur as submodules of irreducible $M^0_k$-modules.
In this subsection we determine the fusion rule of type 
$\binom{L^3}{L^1\ L^2}$.
As a first step, we give a lower bound for every fusion rule by using  Theorem \ref{theorem : main}.
\begin{lemma}\label{lemma : lb}
\begin{enumerate}
\item The fusion rule of following types is greater than or equal to $1$:
Let $i,j\in\{0,1,2\}$. 
\[
\begin{matrix}
\binom{M_k^{0 (i+j)}}{M_k^{0 (i)}\ M_k^{0 (j)}},& \binom{W_k^{0 (i+j)}}{W_k^{0 (i)}\ W_k^{0 (j)}}, & 
\binom{W_k^{0 (i+j)}}{M_k^{0 (i)}\ W_k^{0 (j)}},& \binom{M^a_k}{M_k^{0 (i)}\ M^a_k},\\
\binom{W^a_k}{M_k^{0 (i)}\ W^a_k},& \binom{M^a_k}{W_k^{0 (i)}\ W^a_k},& \binom{W^a_k}{W_k^{0 (i)}\ W^a_k}.
\end{matrix}
\]
\item The fusion rule of following types is greater than or equal to $2$:
\[
\begin{matrix}
\binom{M^a_k}{M^a_k\ M^a_k},& \binom{W^a_k}{M^a_k\ W^a_k},& 
\binom{W^a_k}{W^a_k\ W^a_k}.
\end{matrix}
\]
\end{enumerate} 
\end{lemma}
\begin{proof}
We use notations in Section 3. The fusion rules for $M_k^0$
are obtained in \cite{LY}.
\begin{enumerate}
\item We consider the case that ${\cal S}_1={\cal S}_2={\cal S}_3=\{W_k^0\}$. We have
$\dim_{\C}I_{M_k^0}\binom{W_k^0}{W_k^0\ W_k^0}=1$.
Let $\C x_i\ (i=0,1,2)$ be the one dimensional $G$-modules
such that $\tau\cdot x_i=\xi^ix_i$. Fix nonzero $v^{i}\in \mbox{Hom}_{\C G}(\C x_i, W_k^{0 (i)})$.
Then
\begin{eqnarray*}
W_k^0 & = & 
\bigoplus_{i=0}^{2}x_i\otimes \mbox{Hom}_{\C G}(\C x_i, W_k^{0 (i)}),\\
{\cal I} & = & I_{M_k^0}\binom{W_k^0}{W_k^0\ W_k^0}\otimes_{\C}W_k^0\otimes_{\C}W_k^0,\mbox{ and }\\
{\cal I}_{W_k^{0 (i)},W_k^{0 (j)}}(v^i,v^j) & = &
I_{M_k^0}\binom{W_k^0}{W_k^0\ W_k^0}\otimes (x_i\otimes v^i)\otimes 
(x_j\otimes v^j).
\end{eqnarray*}
Let
$f\in I_{M_k^0}\binom{W_k^0}{W_k^0\ W_k^0}$ be a nonzero intertwining operator. 
Then ${}_{\tau}f=f$ by the construction of intertwining operators 
in \cite{LY}.
So
\[(\tau\otimes e(W_k^0))\cdot f\otimes (x_i\otimes v^i)\otimes 
(x_j\otimes v^j)=\xi^{i+j}f\otimes (x_i\otimes v^i)\otimes 
(x_j\otimes v^j)\]
and $\dim_{\C}I_{\cal W}\binom{W_k^{0 (i+j)}}{W_k^{0 (i)}\ W_k^{0 (j)}}\geq 1$.
We can compute the other cases in the same way.
\item 
We consider the case that ${\cal S}_1={\cal S}_2={\cal S}_3=\{M^a_k,M^b_k,M^c_k\}$.
For $i_1,i_2,i_3\in\{a,b,c\}$, we have
\[
\dim_{\C}I_{M_k^0}\binom{M^{i_3}_k}{M^{i_1}_k\ M^{i_2}_k}
=\left\{\begin{array}{ll}
1, & \mbox{if }\{i_1,i_2,i_3\}=\{a,b,c\},\\
0, & \mbox{otherwise.}
\end{array}\right.
\]
We define an action of $\tau$ on $\{a,b,c\}$ by $\tau(a)=b,\tau(b)=c$ and $\tau(c)=a$.
It is possible to take $\phi(\tau,M^i_k)\ (i=a,b,c)$ such that 
\[\phi(\tau,M^c_k)\phi(\tau,M^b_k)\phi(\tau,M^a_k)=\mbox{id}_{M_k^a}.\]
Fix a nonzero $v^a\in M^a_k$ and 
set $v^b=\phi(\tau, M_k^a)v^a\in M^b_k, v^c=\phi(\tau, M_k^b)v^b\in M^c_k$.
Fix nonzero $f_{a,b,c}\in I_{M_k^0}\binom{M^c_k}{M^a_k\ M^b_k},
f_{b,a,c}\in I_{M_k^0}\binom{M^c_k}{M^b_k\ M^a_k}$
and set 
\[
\begin{array}{rcll}
f_{b,c,a} & = & 
{}_{\tau}f_{a,b,c}\in I_{M_k^0}\binom{M^{a}_k}{M^{b}_k\ M^{c}_k},\\
f_{c,a,b} & = & 
{}_{\tau}f_{b,c,a}\in I_{M_k^0}\binom{M^{b}_k}{M^{c}_k\ M^{a}_k},\\
f_{c,b,a} & = & 
{}_{\tau}f_{b,a,c}\in I_{M_k^0}\binom{M^{a}_k}{M^{c}_k\ M^{b}_k},\\
f_{a,c,b} & = & 
{}_{\tau}f_{c,b,a}\in I_{M_k^0}\binom{M^{b}_k}{M^{a}_k\ M^{c}_k}.
\end{array}
\]
We have
\[{\cal I}=\bigoplus_{\stackrel{\scriptstyle i_1,i_2,i_3\in \{a,b,c\}}{
\{i_1,i_2,i_3\}=\{a,b,c\}}}
I_{M_k^0}\binom{M^{i_3}_k}{M^{i_1}_k\ M^{i_2}_k}
\otimes_{\C}M^{i_1}_k\otimes_{\C}M^{i_2}_k\] 
and
\begin{eqnarray*}
{\cal I}_{M^a_k,M^a_k} & = & \bigoplus_{i=0,1,2}
\C f_{\tau^i(a),\tau^i(b),\tau^i(c)}\otimes v^{\tau^i(a)}\otimes v^{\tau^i(b)}\\
& & {}\oplus\bigoplus_{i=0,1,2}
\C f_{\tau^i(b),\tau^i(a),\tau^i(c)}\otimes v^{\tau^i(b)}\otimes v^{\tau^i(a)}.
\end{eqnarray*}
Since $\bigoplus_{i=0,1,2}
\C f_{\tau^i(a),\tau^i(b),\tau^i(c)}\otimes v^{\tau^i(a)}
\otimes v^{\tau^i(b)}$ and 
$\bigoplus_{i=0,1,2}
\C f_{\tau^i(b),\tau^i(a),\tau^i(c)}\otimes v^{\tau^i(b)}
\otimes v^{\tau^i(a)}$ are isomorphic 
irreducible ${\cal A}_{\alpha_3}(G,{\cal S}_3)$-modules,
we have $\dim_{\C}I_{\cal W}\binom{M^a_k}{M^a_k\ M^a_k}\geq 2$.
We can compute the other cases in the same way.
\end{enumerate}
\end{proof}
We will show that every fusion rule meets the lower bound 
obtained in Lemma \ref{lemma : lb}. We use Lemma \ref{lemma : ub}.
For an irreducible ${\cal W}$-module 
$N=\bigoplus_{n=0}^{\infty}N(n)$, we fix a nonzero vector in $N(0)$
and denote it by $w_{N}$.
By the same argument as in \cite[Lemma 5.2]{DLTYY},
$N$ is spanned by the vectors of 
the form
\begin{eqnarray}\label{eqn : form}
L(-m_1)\cdots L(-m_p)J(-n_1)\cdots J(-n_q)w_{N}
\end{eqnarray}
with $m_1\geq \cdots \geq m_p\geq 1, n_1\geq \cdots \geq n_q\geq 1, p=0,1,\ldots,$
and $q=0,1,\ldots$.
By the same argument as in \cite[Section 5.4]{DLTYY}, 
we have for $n\leq -1$ and $u\in N$,
\begin{eqnarray}\label{eqn : zhu}
L(n)u & = & (-1)^{n-1}[\omega]*[u]+(-1)^{n-1}n[u]*[\omega]+(-1)^n\mbox{\rm wt}u[u],\nonumber\\
J(n)u & = & (-1)^n\Big(n[J(-1)u]+(n+1)[J[0]u]-(n+1)[J]*[u]\nonumber\\
& & {}-\frac{n(n+1)}{2}[u]*[J]\Big)
\end{eqnarray}
in $A(N)$.
Using formula (\ref{eqn : zhu}) and commutation relations (\ref{eqn : com})
repeatedly, it is shown that 
$A(N)$ is generated by $\{J(-1)^iw_{N}\}_{i=0}^{\infty}$ as a $A({\cal W})$-bimodule.

A {\it singular vector} $w$ of weight $h$ for $N$
is by definition a vector $w$ which satisfies
$L(0)w=hw$ and $L(n)w=J(n)w=0$ for $n\geq 1$. 
By commutation relations (\ref{eqn : com}), it is easy to show that 
$w$ is a singular vector
of weight $h$ if and only if $L(0)w=hw$ and $L(1)w=L(2)w=J(1)w=0$. 
If $h-\mbox{wt}N(0)>0$, then
the submodule of $N$ generated by $w$ dose not contain $N(0)$. 
We hence have $w=0$.

\begin{theorem}\label{theorem : fusion}
The fusion rule $L^1\times L^2$ is given by the following list, where 
$\{L^1,L^2\}$ is an arbitrary pair of irreducible ${\cal W}$-modules which
occur as submodules of irreducible $M_k^0$-modules:
Let $i,j\in\{0,1,2\}$.
\begin{eqnarray*}
M_k^{0 (i)}\times M_k^{0 (j)} & = & M_k^{0 (i+j)},\\
M_k^{0 (i)}\times W_k^{0 (j)} & = & W_k^{0 (i+j)},\\
W_k^{0 (i)}\times W_k^{0 (j)} & = & M_k^{0 (i+j)}+W_k^{0 (i+j)},\\
M_k^{0 (i)}\times M_k^a & = & M_k^a,\\
M_k^{0 (i)}\times W_k^a & = & W_k^a,\\
W_k^{0 (i)}\times M_k^a & = & W_k^a,\\
W_k^{0 (i)}\times W_k^a & = & M_k^a+W_k^a,\\
M_k^a\times M_k^a & = & \sum_{i=0}^{2}M_k^{0 (i)}+2M_k^a,\\
M_k^a\times W_k^a & = & \sum_{i=0}^{2}W_k^{0 (i)}+2W_k^a,\\
W_k^a\times W_k^a & = & \sum_{i=0}^{2}M_k^{0 (i)}+\sum_{i=0}^{2}W_k^{0 (i)}+2M_k^a+2W_k^a.
\end{eqnarray*}
\end{theorem}
\begin{proof}
For an irreducible ${\cal W}$-module $N$,
$h_N$ denotes the eigenvalue for $L(0)$ on $N(0)$ and
$k_N$ denotes the eigenvalue for $J(0)$ on $N(0)$.
For all irreducible ${\cal W}$-modules those eigenvalues are computed in \cite{DLTYY}.
For irreducible ${\cal W}$-modules $L^i\ (i=1,2,3)$, set
\begin{eqnarray*}
F(L^1,L^2,L^3) & = & \dim_{\C}\big(L^3(0)^{*}\otimes_{A({\cal W})}A(L^1)
\otimes_{A({\cal W})}L^2(0)\big).
\end{eqnarray*}
Note that $F(L^1,L^2,L^3)$ is lower than or equal to the number of generators
of $A(L^1)$ as a $A({\cal W})$-bimodule since the dimension of 
top level of every irreducible ${\cal W}$-module is $1$.

The simplicity of ${\cal W}$ implies that
\begin{eqnarray*}
\dim_{\C}I_{\cal W}\binom{L^3}{{\cal W}\ L^2} & = & 
\left\{\begin{array}{ll}
1, & \mbox{if }L^2\simeq L^3,\\
0, & \mbox{otherwise}
\end{array}\right.
\end{eqnarray*}
for all irreducible ${\cal W}$-modules $L^2$ and $L^3$.
We consider the case $L^1\neq {\cal W}$.
We will show that $F(L^1,L^2,L^3)$ is lower than or equal to the lower bound
for $\dim_{\C}I_{\cal W}\binom{L^3}{L^1\ L^2}$ given in Lemma \ref{lemma : lb}.
Then, by Lemma \ref{lemma : sym} and Lemma \ref{lemma : ub} we get the desired results.
We use a computer algebra system Risa/Asir to find singular vectors used in the following
argument.The explicit forms of those singular vectors are given in Appendix. 
\begin{enumerate}
\item
We consider the case $L^1=W_k^{0 (1)}$.
Since $(5\sqrt{-3}L(-1)+J(-1))w_{W_k^{0 (1)}}$ is a singular vector and
$A(W_k^{0 (1)})$ is generated by 
$\{J(-1)^iw_{W_k^{0 (1)}}\}_{i=0}^{\infty}$ as a $A({\cal W})$-bimodule,
$A(W_k^{0 (1)})$ is generated by $w_{W_k^{0 (1)}}$ as a $A({\cal W})$-bimodule.
So $F(W_k^{0 (1)},L^2,L^3)\leq 1$.
Set 
\begin{eqnarray*}
w_1 & = & \big(-30\sqrt{-3}L(-1)J(-1)+39\sqrt{-3}J(-2)\\
& & {}+5J(-1)^2+336L(-2)+405L(-1)^2\big)w_{W_k^{0 (1)}}.
\end{eqnarray*}
Since $w_1$ is a singular vector, we have a relation
\begin{eqnarray*}
0 & = & 
50\sqrt{-3}([\omega]^2*[w_{W_k^{0 (1)}}]+[w_{W_k^{0 (1)}}]*[\omega]^2)\\
& & -20\sqrt{-3}([\omega]*[w_{W_k^{0 (1)}}]+[w_{W_k^{0 (1)}}]*[\omega])\\
& & {}+4\sqrt{-3}[w_{W_k^{0 (1)}}]-100\sqrt{-3}[\omega]*[w_{W_k^{0 (1)}}]*[\omega]\\
& & {}-5[w_{W_k^{0 (1)}}]*[J]+5[J]*[w_{W_k^{0 (1)}}]
\end{eqnarray*}
in $A(W_0^k(1))$ by using  formulas (\ref{eqn : zhu}). 
Therefore,
\begin{eqnarray*}
0 & = & 
\big(50\sqrt{-3}(h_{L^2}^2+h_{L^3}^2)-20\sqrt{-3}(h_{L^2}+h_{L^3})\\
& & {}+4\sqrt{-3}-100\sqrt{-3}h_{L^2}h_{L^3}-5k_{L^2}+5k_{L^3}\big)
w_{(L^3)^{\prime}}\otimes [w_{W_k^{0 (1)}}]\otimes w_{L^2}
\end{eqnarray*}
in 
$L^3(0)^{*}\otimes_{A({\cal W})}A(W_k^{0 (1)})\otimes_{A({\cal W})}L^2(0)$.
Set 
\begin{eqnarray*}
\psi(h_{L^2},k_{L^2},h_{L^2},k_{L^3}) & = & 
50\sqrt{-3}(h_{L^2}^2+h_{L^2}^2)-20\sqrt{-3}(h_{L^2}+h_{L^3})\\
& & {}+4\sqrt{-3}-100\sqrt{-3}h_{L^2}h_{L^3}-5k_{L^2}+5k_{L^3}.
\end{eqnarray*}
If $\psi(h_{L^2},k_{L^2},h_{L^2},k_{L^3})\neq 0$, then 
$w_{(L^3)^{\prime}}\otimes [w_{W_k^{0 (1)}}]\otimes w_{L^2}=0$ and 
$F(W_k^{0 (1)},L^2, L^3)=0$. By computing 
$\psi(h_{L^2},k_{L^2},h_{L^3},k_{L^3})$ for all pairs $(L^2, L^3)$ of 
$20$ irreducible ${\cal W}$-modules, we have $F(W_k^{0 (1)},L^2,L^3)\leq 1$ if the pair
$(L^2,L^3)$ is one of 
\[
\begin{array}{ll}
(M_k^{0 (i)},W_k^{0 (1+i)}),(W_k^{0 (i)},M_k^{0 (1+i)}), 
(W_k^{0 (i)},W_k^{0 (1+i)}), & (i=0,1,2),\\
(M_k^a,W_k^a),(W_k^a,M_k^a), \mbox{ and }(W_k^a,W_k^a)
\end{array}\]
and $F(W_k^{0 (1)},L^2,L^3)=0$ otherwise. Combining these results and 
Lemma \ref{lemma : lb}, we have 
\begin{eqnarray*}
W_k^{0 (1)}\times M_k^{0 (i)} & = & W_k^{0 (1+i)},\\
W_k^{0 (1)}\times W_k^{0 (i)} & = & M_k^{0 (1+i)}+W_k^{0 (1+i)},\\
W_k^{0 (1)}\times M_k^a & = & W_k^a,\\
W_k^{0 (1)}\times W_k^a & = & M_k^a+W_k^a.
\end{eqnarray*}

In the case $L^1=W_k^{0 (2)}$, we can compute the fusion rules
in the same way.
In the case $L^1=M_k^{0 (i)}\ (i=1,2)$, 
it is shown that $A(M_k^{0 (i)})$ is generated by $[w_{M_k^{0 (i)}}]$
in the same way. But in these cases we need two
singular vectors $v_{31}\in M_k^{0 (i)}(3)$ and $v_{41}\in M_k^{0 (i)}(4)$.
We can compute the fusion rules by using  two relations 
$v_{41}=0$ and $J(-1)v_{31}=0$.
\item
We consider the case $L^1=W_k^{0 (0)}$.
There are two singular vectors $u_{2i}\in W_k^{0 (0)}(2)\ (i=1,2)$ 
and there is one singular vector $u_{41}\in W_k^{0 (0)}(4)$.
Since $u_{21}=0$, $A(W_0^k(0))$ is generated by 
$\{w_{W_k^{0 (0)}}, J(-1)w_{W_k^{0 (0)}}\}$ as a $A({\cal W})$-bimodule
and for $X\in W_k^{0 (0)}$ we have 
$a_{1}(L^2,L^3;X), a_2(L^2,L^3;X)\in\C$ such that
\begin{eqnarray*}
w_{(L^3)^{\prime}}\otimes 
[X]\otimes w_{L^2} & = & a_{1}(L^2,L^3;X)
w_{(L^3)^{\prime}}\otimes 
[w_{W_k^{0 (0)}}]\otimes w_{L^2}\\
& & {}+a_{2}(L^2,L^3;X)
w_{(L^3)^{\prime}}\otimes 
[J(-1)w_{W_k^{0 (0)}}]\otimes w_{L^2}.
\end{eqnarray*}
Therefore, $F(W_k^{0 (0)}, L^2, L^3)\leq 2$.
Set a matrix 
\begin{eqnarray*}
A & = & \begin{pmatrix}
a_{1}(L^2,L^3;u_{22}) & a_{2}(L^2,L^3;u_{22})\\
a_{1}(L^2,L^3;J(-1)u_{22}) & a_{2}(L^2,L^3;J(-1)u_{22})\\
a_{1}(L^2,L^3;v_{41}) & a_{2}(L^2,L^3;u_{41})\\
a_{1}(L^2,L^3;J(-1)u_{41}) & a_{2}(L^2,L^3;J(-1)u_{41})
\end{pmatrix}.
\end{eqnarray*}
We have 
\[
A\begin{pmatrix}
(w_{(L^3)^{\prime}}\otimes[w_{W_k^{0 (0)}}]\otimes w_{L^2}\\
w_{(L^3)^{\prime}}\otimes[J(-1)w_{W_k^{0 (0)}}]\otimes w_{L^2}\end{pmatrix}
=\begin{pmatrix}0\\0\\0\\0\end{pmatrix}.
\]
If $\mbox{rank}A=2$, then 
$w_{(L^3)^{\prime}}\otimes 
[w_{W_k^{0 (0)}}]\otimes w_{L^2}=w_{(L^3)^{\prime}}\otimes 
[J(-1)w_{W_k^{0 (0)}}]\otimes w_{L^2}=0$ in 
$L^3(0)^{*}\otimes_{A({\cal W})}A(W_k^{0 (0)})\otimes_{A({\cal W})}L^2(0)$.
If $\mbox{rank}A=1$, then 
$w_{(L^3)^{\prime}}\otimes 
[w_{W_k^{0 (0)}}]\otimes w_{L^2}$ and $
w_{(L^3)^{\prime}}\otimes [J(-1)w_{W_k^{0 (0)}}]\otimes w_{L^2}$ 
are linearly dependent.
We hence have 
\begin{eqnarray*}
F(W_k^{0 (0)}, L^2, L^3) & \leq &
\left\{ 
\begin{array}{ll}
0, & \mbox{if }\mbox{rank}A=2,\\
1, & \mbox{if }\mbox{rank}A=1,\\
2, & \mbox{if }\mbox{rank}A=0.
\end{array}\right.
\end{eqnarray*}
By computing the rank of $A$ for all pairs $(L^2, L^3)$ of 
$20$ irreducible ${\cal W}$-modules, we have for $i=0,1,2$
\begin{eqnarray*}
W_k^{0 (0)}\times M_k^{0 (i)} & = & W_k^{0 (i)},\\
W_k^{0 (0)}\times W_k^{0 (i)} & = & M_k^{0 (i)}+W_k^{0 (i)},\\
W_k^{0 (0)}\times M_k^a & = & W_k^a,\\
W_k^{0 (0)}\times W_k^a & = & M_k^a+W_k^a.
\end{eqnarray*}
In the case that $L^1$ is one of $M_k^a$ and $W_k^a$,
we can compute the fusion rules is the same way. 
We roughly explain each case. 
In the case $L^1=M_k^a$,
there are two singular vectors $v_{21}, v_{22}$ in $M_k^a(2)$ and 
three singular vectors $v_{61},v_{62},v_{63}$ in $M_k^a(6)$.
Set the same matrix for $v_{22}, v_{61},v_{62},v_{63},
J(-1)v_{61},J(-1)v_{62}$ and $J(-1)v_{63}$ as in the case of $W_k^{0 (0)}$.
By computing the rank of the matrix, we can determine the fusion rules.
In the case $L^1=W_k^a$, there is a singular vector $v_{21}$ in $W_k^a(2)$ and 
there are two singular vectors $v_{41},v_{42}$ in $W_k^a(4)$.
Set the same matrix for $v_{41},v_{42},J(-1)v_{41}$ and $J(-1)v_{42}$ as 
in the case of $W_k^{0 (0)}$.
By computing the rank of the matrix, we can determine the fusion rules.
\end{enumerate}
\end{proof}
\section{Appendix}
We give some singular vectors in irreducible ${\cal W}$-modules used in 
Theorem \ref{theorem : fusion}.
We omit $w_N\in N(0)$ from the explicit form of every singular vector in a
${\cal W}$-module $N$. 

\begin{itemize}
\item $M_k^{0 (1)}$.
\[
\begin{array}{cl}
(1) & 9\sqrt{-3}L(-1)+J(-1),\\
(2) & 6255L(-1)J(-2)-375\sqrt{-3}L(-1)J(-1)^2+36960L(-2)J(-1)\\
& {}+5175L(-1)^2J(-1)-26208J(-3)+1425\sqrt{-3}J(-2)J(-1)\\
& {}+25J(-1)^3+23040\sqrt{-3}L(-3)+147600\sqrt{-3}L(-2)L(-1)\\
& {}-44625\sqrt{-3}L(-1)^3,\\
(3) & 
-36720\sqrt{-3}L(-1)J(-3)+324L(-1)J(-2)J(-1)\\
& {}-16\sqrt{-3}L(-1)J(-1)^3+9360\sqrt{-3}L(-2)J(-2)\\
& {}+12852\sqrt{-3}L(-1)^2J(-2)+2400L(-2)J(-1)^2\\
& {}+462L(-1)^2J(-1)^2-5040\sqrt{-3}L(-3)J(-1)\\
& {}+8640\sqrt{-3}L(-2)L(-1)J(-1)-5232\sqrt{-3}L(-1)^3J(-1)\\
& {}+35280\sqrt{-3}J(-4)-819J(-2)^2\\
& {}+76\sqrt{-3}J(-2)J(-1)^2+J(-1)^4\\
& {}-751680)L(-4)-1028160)L(-3)L(-1)\\
& {}+254880L(-2)L(-1)^2+16929L(-1)^4.
\end{array}
\]
We obtain singular vectors in $M_k^{0 (2)}$ by replacing 
$J(n)$ with $-J(n)$ in the above vectors. 
\item
$W_k^{0 (0)}$.
\[
\begin{array}{cl}
(1) & -70\sqrt{-3}L(-1)J(-1)+91\sqrt{-3}J(-2)-5J(-1)^2-2496L(-2)+195L(-1)^2,\\
(2) & -70\sqrt{-3}L(-1)J(-1)+91\sqrt{-3}J(-2)+5J(-1)^2+2496L(-2)-195L(-1)^2,\\
(3) & 
-1500L(-1)J(-2)J(-1)+1200L(-2)J(-1)^2+750L(-1)^2J(-1)^2\\
& {}+3600J(-3)J(-1)+825J(-2)^2+J(-1)^4\\
& {}-633600L(-4)+46800L(-3)L(-1)+230400L(-2)^2\\
& {}-126000L(-2)L(-1)^2+50625L(-1)^4,
\end{array}
\]
\item 
$W_k^{0 (1)}$.
\[
\begin{array}{cl}
(1) & 5\sqrt{-3}L(-1)+J(-1),\\
(2) & -30\sqrt{-3}L(-1)J(-1)+39\sqrt{-3}J(-2)+5J(-1)^2+336L(-2)
+405L(-1)^2.
\end{array}
\]
We obtain singular vectors in $W_k^{0 (2)}$ by replacing 
$J(n)$ with $-J(n)$ in the above vectors. 
\item
$W_k^a$.
\[
\begin{array}{cl}
(1) & J(-1)^2-30L(-2)+75L(-1)^2,\\
(2) & -5040\sqrt{-3}L(-1)J(-3)-4980L(-1)J(-2)J(-1)\\
& {}+200\sqrt{-3}L(-1)J(-1)^3+7668\sqrt{-3}L(-2)J(-2)\\
& {}+990\sqrt{-3}L(-1)^2J(-2)+9660L(-2)J(-1)^2\\
& {}+7350L(-1)^2J(-1)^2+35760\sqrt{-3}L(-3)J(-1)\\
& {}-78960\sqrt{-3}L(-2)L(-1)J(-1)+34200\sqrt{-3}L(-1)^3J(-1)\\
& {}+5208\sqrt{-3}J(-4)+12780J(-3)J(-1)\\
& {}+2622J(-2)J(-2)-310\sqrt{-3}J(-2)J(-1)^2\\
& {}+25J(-1)^4-9000L(-4)\\
& {}-255780L(-3)L(-1)-28044L(-2)^2\\
& {}+557460L(-2)L(-1)^2-78975L(-1)^4,\\
(3) & 5040\sqrt{-3}L(-1)J(-3)-4980L(-1)J(-2)J(-1)\\
& {}-200\sqrt{-3}L(-1)J(-1)^3-7668\sqrt{-3}L(-2)J(-2)\\
& {}-990\sqrt{-3}L(-1)^2J(-2)+9660L(-2)J(-1)^2\\
& {}+7350L(-1)^2J(-1)^2-35760\sqrt{-3}L(-3)J(-1)\\
& {}+78960\sqrt{-3}L(-2)L(-1)J(-1)-34200\sqrt{-3}L(-1)^3J(-1)\\
& {}-5208\sqrt{-3}J(-4)+12780J(-3)J(-1)\\
& {}+2622J(-2)J(-2)+310\sqrt{-3}J(-2)J(-1)^2\\
& {}+25J(-1)^4-9000L(-4)\\
& {}-255780L(-3)L(-1)-28044L(-2)^2\\
& {}+557460L(-2)L(-1)^2-78975L(-1)^4.
\end{array}
\]
\item $M_k^a$.
\[
\begin{array}{cl}
(1) & 8\sqrt{-3}L(-1)J(-1)-6\sqrt{-3}J(-2)+J(-1)^2+90L(-2)+27L(-1)^2,\\
(2) & -8\sqrt{-3}L(-1)J(-1)+6\sqrt{-3}J(-2)+J(-1)^2+90L(-2)+27L(-1)^2,
\end{array}
\]
There are another three singular vectors $v_{6i}\ (i=1,2,3)$ in $M_k^a$,

\[
\begin{array}{cl}
v_{61}=& 63987840\sqrt{-3}L(-1)J(-5)-1059480L(-1)J(-4)J(-1)\\ & {}
+1587600L(-1)J(-3)J(-2)+14400\sqrt{-3}L(-1)J(-3)J(-1)^2\\ & {}
-26208\sqrt{-3}L(-1)J(-2)^2J(-1)+4260L(-1)J(-2)J(-1)^3\\ & {}
+32\sqrt{-3}L(-1)J(-1)^5+23284800\sqrt{-3}L(-2)J(-4)\\ & {}
-19867680\sqrt{-3}L(-1)^2J(-4)+145800L(-2)J(-3)J(-1)\\ & {}
-1414260L(-1)^2J(-3)J(-1)-427140L(-2)J(-2)^2\\ & {}
-155358L(-1)^2J(-2)^2-15840\sqrt{-3}L(-2)J(-2)J(-1)^2\\ & {}
+10224\sqrt{-3}L(-1)^2J(-2)J(-1)^2+3990L(-2)J(-1)^4\\ & {}
-447L(-1)^2J(-1)^4-25401600\sqrt{-3}L(-3)J(-3)\\ & {}
-5443200\sqrt{-3}L(-2)L(-1)J(-3)+1995840\sqrt{-3}L(-1)^3J(-3)\\ & {}
-659880L(-3)J(-2)J(-1)+1666440L(-2)L(-1)J(-2)J(-1)\\ & {}
+341388L(-1)^3J(-2)J(-1)-41280\sqrt{-3}L(-3)J(-1)^3\\ & {}
+32640\sqrt{-3}L(-2)L(-1)J(-1)^3+7872\sqrt{-3}L(-1)^3J(-1)^3\\ & {}
+4294080\sqrt{-3}L(-4)J(-2)+16420320\sqrt{-3}L(-3)L(-1)J(-2)\\ & {}
+4989600\sqrt{-3}L(-2)^2J(-2)-1360800\sqrt{-3}L(-2)L(-1)^2J(-2)\\ & {}
+122472\sqrt{-3}L(-1)^4J(-2)-3483720L(-4)J(-1)^2\\ & {}
-1742220L(-3)L(-1)J(-1)^2+38700L(-2)^2J(-1)^2\\ & {}
+588420L(-2)L(-1)^2J(-1)^2-169749L(-1)^4J(-1)^2\\ & {}
-42370560\sqrt{-3}L(-5)J(-1)-45861120\sqrt{-3}L(-4)L(-1)J(-1)\\ & {}
-2332800\sqrt{-3}L(-3)L(-2)J(-1)-31164480\sqrt{-3}L(-3)L(-1)^2J(-1)\\ & {}
-3542400\sqrt{-3}L(-2)^2L(-1)J(-1)+4429440\sqrt{-3}L(-2)L(-1)^3J(-1)\\ & {}
+23328\sqrt{-3}L(-1)^5J(-1)-101787840\sqrt{-3}J(-6)\\ & {}
-546480J(-5)J(-1)-1186920J(-4)J(-2)\\ & {}
+5280\sqrt{-3}J(-4)J(-1)^2-2381400J(-3)^2\\ & {}
-47520\sqrt{-3}J(-3)J(-2)J(-1)-3420J(-3)J(-1)^3\\ & {}
+11088\sqrt{-3}J(-2)^3+870J(-2)^2J(-1)^2\\ & {}
-88\sqrt{-3}J(-2)J(-1)^4+J(-1)^6\\ & {}
+879076800L(-6)+990072720L(-5)L(-1)\\ & {}
+65091600L(-4)L(-2)+323666280L(-4)L(-1)^2\\ & {}
+102173400L(-3)^2+73823400L(-3)L(-2)L(-1)\\ & {}
-3027780L(-3)L(-1)^3-152523000L(-2)^3\\ & {}
+95652900L(-2)^2L(-1)^2-15326010L(-2)L(-1)^4\\ & {}
-505197L(-1)^6.
\end{array}
\]
The vector $v_{62}$ is obtained by replacing 
$J(n)$ with $-J(n)$ in $v_{61}$.
\[
\begin{array}{ll}
v_{63}=& 
-1478136600L(-1)J(-4)J(-1)+423979920L(-1)J(-3)J(-2)\\ & {}
-1273500L(-1)J(-2)J(-1)^3-29005560L(-2)J(-3)J(-1)\\ & {}
+58538700L(-1)^2J(-3)J(-1)+134322300L(-2)J(-2)^2\\ & {}
-133840350L(-1)^2J(-2)^2+1341750L(-2)J(-1)^4\\ & {}
+680625L(-1)^2J(-1)^4-778588200L(-3)J(-2)J(-1)\\ & {}
+143310600L(-2)L(-1)J(-2)J(-1)+37975500L(-1)^3J(-2)J(-1)\\ & {}
-2232505800L(-4)J(-1)^2-311863500L(-3)L(-1)J(-1)^2\\ & {}
+586133100L(-2)^2J(-1)^2-239341500L(-2)L(-1)^2J(-1)^2\\ & {}
+127186875L(-1)^4J(-1)^2-393666480J(-5)J(-1)\\ & {}
-1236672360J(-4)J(-2)-104786136J(-3)^2\\ & {}
+3968100J(-3)J(-1)^3+2139750J(-2)^2J(-1)^2\\ & {}
+625J(-1)^6+275225065920L(-6)\\ & {}
+450006898320L(-5)L(-1)+221829042960L(-4)L(-2)\\ & {}
-60223224600L(-4)L(-1)^2+121440173400L(-3)^2\\ & {}
+145205865000L(-3)L(-2)L(-1)-22975690500L(-3)L(-1)^3\\ & {}
-5826151800L(-2)^3+65516566500L(-2)^2L(-1)^2\\ & {}
-55664516250L(-2)L(-1)^4+5974171875L(-1)^6.
\end{array}
\]
\end{itemize}


\begin{thebibliography}{99}
\bibitem{B}
R. Borcherds, 
Vertex algebras, Kac-Moody algebras, and the Monster, 
{\em Proc. Nat. Acad. Sci. U.S.A.} {\bf 83} (1986), 3068--3071.
\bibitem{DVVV}
R. Dijkgraaf, C. Vafa, E. Verlinde, and H. Verlinde, 
The operator algebra of orbifold models, 
{\em Comm. Math. Phys.} {\bf 123} (1989), 485--526.
\bibitem{DLTYY}
 C. Dong, C.H. Lam, K. Tanabe, H. Yamada, and K. Yokoyama,
$\Z_3$ symmetry and $W_3$ algebra in lattice
vertex operator algebras, 
{\em Pacific J. Math.} {\bf 215} (2004), 245--296.
\bibitem{DL}
C. Dong and J. Lepowsky, Generalized Vertex Algebras and Relative Vertex Operators,
Progress in Math., Vol. {\bf 112}, Birkhauser, Boston, 1993.
\bibitem{DM}
C. Dong and G. Mason, 
On quantum Galois theory, 
{\em Duke Math. J.} {\bf 86} (1997), 305--321.
\bibitem{DY}
C. Dong and G. Yamskulna,
Vertex operator algebras, generalized doubles and dual pairs,
{\em Math. Z.} {\bf 241} (2002), 397--423.
\bibitem{FaZ}
V. A. Fateev and A. B. Zamolodchikov, Conformal quantum field
theory models in two dimensions having $\Z_3$ symmetry, {\itshape
Nuclear Physics} {\bfseries B280} (1987), 644--660.
\bibitem{FHL}
I. Frenkel, Y. Huang, and J. Lepowsky, 
On axiomatic approaches to vertex operator algebras 
and modules, Mem. Amer. Math. Soc. {\bf 104} (1993).
\bibitem{FLM}
I. B. Frenkel, J. Lepowsky and A. Meurman, {\itshape Vertex
Operator Algebras and the Monster}, Pure and Applied Math., Vol.
{\bfseries 134}, Academic Press, 1988.
\bibitem{FZ}
I.\ Frenkel and Y.\ Zhu, 
Vertex operator algebras associated to representations of
affine and Virasoro algebras, {\em Duke Math. J.} {\bf 66} (1992), 123--168.
\bibitem{KLY2}
K. Kitazume, C. Lam and H. Yamada, $3$-state Potts model,
moonshine vertex operator algebra and $3A$ elements of the monster
group, {\itshape IMRS}, {\bfseries 23} (2003), 1269--1303.
\bibitem{LY}
C. Lam and H. Yamada, $\mathbb{Z}_{2}\times
\mathbb{Z}_{2}$ codes and vertex operator algebras, {\itshape J.
Algebra} {\bfseries 224} (2000), 268--291.
\bibitem{L}
H.\ Li, 
Determining fusion rules by $A(V)$-modules and bimodules,
{\it J. Algebra} {\bf 212} (1999), 515--556. 
\bibitem{Z}
Y. Zhu, 
Modular invariance of characters of vertex operator algebras,
{\em J. Amer. Math. Soc.} {\bf9} (1996), 237--302. 
\end{thebibliography}
\end{document}